\documentclass[12pt,a4wide]{article}
\usepackage{amsmath, amsfonts, amssymb, amsthm, euscript, amscd, latexsym}
 \usepackage[pdftex,colorlinks=true,
                      pdfstartview=FitV,
                      linkcolor=blue,
                      citecolor=blue,
                      urlcolor=blue,
          ]{hyperref}
 \usepackage[dvips]{color}
\def\eee#1{ \begin{equation} #1 \end{equation} }
\def\aa#1{ \begin{align*} #1 \end{align*} }
\def\aaa#1{ \begin{align} #1 \end{align} }
\def\mm#1{ \begin{multline*} #1 \end{multline*} }
\def\mmm#1{ \begin{multline} #1 \end{multline} }

\newtheorem{thm}{\sc Theorem}
\newtheorem{lem}{\sc Lemma}

\newtheorem{pro}{\sc Proposition}

\newcommand{\sss}{\scriptscriptstyle}

\newcommand{\eps}{\varepsilon}
\newcommand{\pl}{\partial}
\newcommand{\gt}{\geqslant}
\newcommand{\lt}{\leqslant}
\newcommand{\sub}{\subset}

\newcommand{\dl}{\delta}
\newcommand{\al}{\alpha}

 \newcommand{\Dl}{\Delta}
 \newcommand{\la}{\lambda}

\newcommand{\om}{\omega}
\newcommand{\mc}{\mathcal}
\newcommand{\Om}{\Omega}

\newcommand{\td}{\tilde}

\newcommand{\e}{{\sss E}}

\newcommand{\x}{\times}
\newcommand{\mto}{\mapsto}

\newcommand{\W}{\mathbb W}

\newcommand{\C}{{\rm C}}

\DeclareMathOperator{\ind}{\mathbb I}
\newcommand{\lap}{\Delta}
\newcommand{\nab}{\nabla}
\newcommand{\fdot}{\,\cdot\,}
\def\Rnu{{\mathbb R}}

\def\ffi{\varphi}

\def\suml{\sum\limits}
\def\intl{\int\limits}

\def\com#1{}

\long\def\symbolfootnote[#1]#2{\begingroup%
\def\thefootnote{\fnsymbol{footnote}}\footnote[#1]{#2}\endgroup}

\include{srctex.sty}

\begin{document}

\author{Evelina Shamarova
\footnote{Institute for mathematical Methods in Economics, Vienna University of Technology}
 }

\title{\textbf{Chernoff's theorem \\ for evolution families}}

 \maketitle

 \vspace{-7mm}

 \begin{abstract}
 A generalized version of Chernoff's theorem
 has been obtained.
 Namely, the version of Chernoff's theorem for semigroups obtained in
 a paper by Smolyanov, Weizs\"acker, and Wittich \cite{sm_vr} is
 generalized for a time-inhomogeneous case.
 The main theorem obtained in the current paper,
 Chernoff's theorem for evolution families, deals
 with a family of time-dependent generators of semigroups $A_t$
 on a Banach space,
 a two-parameter family of operators $Q_{t,t+\Dl t}$
 satisfying the relation:
 $\left.\frac{\pl}{\pl \Dl t}Q_{t,t+\Dl t}\right|_{\Dl t =
 0}=A_t$, whose products
 $Q_{t_i,t_{i+1}}\dots Q_{t_{k-1},t_k}$
 are uniformly bounded for all subpartitions
 $s = t_0 < t_1 < \cdots < t_n = t$.
 The theorem states that $Q_{t_0,t_1}\dots
 Q_{t_{n-1},t_n}$ converges
 to an evolution family $U(s,t)$ solving a non-autonomous Cauchy problem.
 Furthermore, the theorem is formulated for a particular case when the generators
 $A_t$ are time dependent second order differential operators.
 Finally, an example of application of this theorem to a construction of time-inhomogeneous
 diffusions on a compact Riemannian manifold is given.

\vspace{2mm}

 \noindent \textit{Keywords:} Chernoff's theorem,
 evolution family, strongly continuous semigroup,
 evolution families generated by manifold valued stochastic processes.
 \end{abstract}

  \symbolfootnote[0]{This work was supported by the research grant
  of the Erwin Schr\"odinger Institute for mathematical physics,
  and by the Austrian Science Fund (FWF) under START-prize-grant Y328.
  }
   \symbolfootnote[0]{{Email:\color{blue}\textit{
  Evelina.Shamarova@fam.tuwien.ac.at}}}

 \section{Chernoff's theorem for evolution families}

  \subsection{Notation}
  Let
  $A_t$, $t\in [S,T] \sub \Rnu_+ \cup \{0\}$,
  be generators of strongly continuous semigroups on
  a Banach space $E$. Let $D(A_t)$ denote the domain of $A_t$.
  We assume that there exists
  a Banach space $Y\sub \cap_{t\in [S,T]} D(A_t)$,
  which is dense in $E$.

  Given a $t\in (S,T]$, and an $x\in Y$, we
  consider a non-autonomous Cauchy problem
  on the interval $[S,t]$
  with the final condition $x$:
 \eee{
 \label{Cauchy2}
 \left\{
 \begin{aligned}
 &\dot u(s)=-A_s u(s)\\
 & u(t)=x.
 \end{aligned}
 \right.
 }
 Let $D_{[S,T]}= \{(s,t): \,  s\lt t, \, s\in [S,T], \, t\in [S,T]\}$.
 The evolution family $U(s,t)$, $(s,t)\in D_{[S,T]}$,
 solving the Cauchy problem \eqref{Cauchy2} satisfies the
 relation:
 \eee{
 \label{evol1}
 U(s,r)U(r,t)=U(s,t)
 }
 for all $s\lt r\lt t$ (see \cite{Dalecky-Fomin}, Chapter VI, paragraph 2).
 \subsection{Case of non-commuting generators}
  \renewcommand{\labelenumi}{\theenumi)}

 \begin{thm}[Chernoff's theorem for evolution families]
 \label{chernnon}
 Let $A_t$ be generators of strongly continuous semigroups,
 $Q_{t_1,t_2}$, $t_1,t_2>0$, be a
 two-parameter family of bounded operators $E\to E$,
 and $U(s,t)$, $S\lt s \lt t \lt T$, be an evolution family
 of operators
 with the generators $A_t$ (see \cite{Dalecky-Fomin}, Chapter VI, paragraph 2).
 We assume that the following assumptions are fulfilled:
 \begin{enumerate}
 \item
 \label{well-posedness}
 there exists a Banach space $Y\sub \cap_{t\in [S,T]}D(A_t)$
 which is dense in $E$, invariant under the action of $U(s,t)$ for all
 $(s,t)\in D_{[S,T]}$, i.e.
 $U(s,t)Y\sub Y$, and
 such that
 the Cauchy problem \eqref{Cauchy2}
 is well-posed (backward solvable) for all $x\in Y$;
 \item
 \label{continuous}
 the function
 $[S,t]\to E \,,\, s \mto \frac{\pl}{\pl s}U(s,t)x$
 is continuous for all $x\in Y$ and $t\in [S,T]$;
 \item
 \label{Qcond}
 for any subinterval $[s,t]\sub [S,T]$,
 there exists a constant
 $M(s,t)>0$ such that
 for all sequences
 $\{s = \tau_1 < \tau_2 < \cdots < \tau_ k \lt t \}$,
 $\|Q_{\tau_1,\tau_2}\cdots Q_{\tau_{k-1},\tau_k}\|\lt M(s,t)$;\,\footnote{For
 example, this assumption is fulfilled when $Q_{\tau,\tau+\Dl\tau}$ are
 contractions.}
 \item
 \label{ch_ass}
 for any subinterval $[s,t]\in (S,T]$,
 for any fixed $x\in Y$,
 \eee{
 \label{ass_convergence}
 \frac{Q_{\tau-\Dl\tau,\tau}-I}{\Dl\tau}\, U(\tau,t)x
 \to
 A_{\tau}\, U(\tau,t)x, \quad   \Dl\tau\to 0
 }
 uniformly in $\tau\in [s,t]$.
 \end{enumerate}
 Then, for any subinterval $[s,t]\sub [S,T]$, for any sequence of
 partitions
 $\{s=t_0 < t_1 < \cdots < t_n=t\}$
 of $[s,t]$ such that
 $\max \,(t_{j+1}-t_j)\to 0$ as $n\to\infty$,
 and for all $x\in E$,
\[
 Q_{t_0,t_1}\ldots Q_{t_{n-1},t_n}x
 \to U(s,t)\,x, \quad n\to \infty.
 \]
 \end{thm}
 \begin{proof}
 First we consider the case $s>S$, i.e. $[s,t]\sub (S,T]$.
 We fix an arbitrary $x\in Y$.
 Using relation \eqref{evol1}, we obtain:
 \mmm{
 \label{4}
 Q_{t_0,t_1}Q_{t_1,t_2}\ldots Q_{t_{n-1},t_n}-
 U(s,t) \\ =
 \suml_{j=0}^{n-1}Q_{t_1,t_2}\ldots Q_{t_{j-2},t_{j-1}}
 (Q_{t_{j-1},t_{j}}-U(t_{j-1},t_j))
 U(t_j,t).
 }
 Let $\dl_n=\max_j(t_j-t_{j-1})$
 be the mesh of the partition $\{s=t_0 < t_1 < \cdots < t_n=t\}$.
 Relation \eqref{4} implies the following inequality:
 \aaa{
 \notag
 &\quad\|(Q_{t_0,t_1}Q_{t_1,t_2}\ldots Q_{t_{n-1},t_n}-
 U(s,t))x\|\\
 \notag
 &\lt
 \sum_{j=0}^{n-1}\Dl t_j
 \Bigl\|
 \Bigl(
 \frac{Q_{t_{j-1},t_j}-I}{t_j-t_{j-1}}
 -\frac{U(t_{j-1},t_j)-I}{t_j-t_{j-1}}
 \Bigr)U(t_j,t)\,x
 \Bigr\|_{\sss E}\\
 \notag
 &\lt
 (t-s)\sup
 \Bigl\{
 \Bigl\|
 \Bigl(
 \frac{Q_{\tau-\Dl\tau,\tau}-I}{\Dl \tau}
 -\frac{U(\tau-\Dl\tau,\tau)-I}{\Dl \tau}
 \Bigr)
 U(\tau,t)\, x
 \Bigr\|_{\sss E}:\\
 \notag
 &\phantom{(t-s)\sup
 \Bigl\{
 \Bigl\|
 \Bigl(
 \frac{Q_{\tau-\Dl\tau,\tau}-I}{\Dl \tau}
 -\frac{U(\tau-\Dl\tau,\tau)-I}{\Dl \tau}
 \hspace{4mm}
 }
 \tau \in [s,t], 0<\Dl\tau<\dl_n \Bigr\}\\
 \label{1term}
 &\lt
 (t-s)\sup
 \Bigl\{
 \Bigl\|
 \Bigl(
 \frac{Q_{\tau-\Dl\tau,\tau}-I}{\Dl \tau}
 -A_\tau
 \Bigr)
  U(\tau,t)\,x
 \Bigr\|_{\sss E}:
 \tau \in [s,t], 0<\Dl\tau< \dl_n
 \Bigr\}\\
 \notag
 &+(t-s)\sup
 \Bigl\{
 \Bigl\|
 \frac{U(\tau-\Dl\tau,t)\,x-U(\tau,t)\,x}{\Dl \tau}-A_\tau
 U(\tau,t)\,x
 \Bigr\|_{\sss E}:\\
 \label{3term}
&\phantom{(t-s)\sup
 \Bigl\{
 \Bigl\|
 \Bigl(
 \frac{Q_{\tau-\Dl\tau,\tau}-I}{\Dl \tau}
 -\frac{U(\tau-\Dl\tau,\tau)-I}{\Dl \tau}
 \hspace{4mm}
 }
 \tau \in [s,t), 0<\Dl\tau<\dl_n
 \Bigr\}
 }
 Assumption \ref{ch_ass} implies that
 the term \eqref{1term} converges to zero.
 Let us consider the term \eqref{3term}. Assumption
 \ref{continuous} implies that the function
 \[
 [s,t]\to E, \; \zeta \mto A_\zeta U(\zeta,t)
 \]
 is continuous since $\frac{\pl}{\pl \zeta}U(\zeta,t)=-A_\zeta U(\zeta,t)$.
 Taking into consideration this, we obtain that
 for $\Dl\tau\in (0,s-S)$
 there exists a $\theta\in(0,1)$ such that
 \eee{
 U(\tau-\Dl\tau,t)x=U(\tau,t)x+\Dl\tau
 A_{\tau-\theta\Dl\tau}U(\tau-\theta\Dl\tau,t)x.
 }
 Hence,
 \mm{
 \frac{U(\tau-\Dl\tau,t)\,x-U(\tau,t)\,x}{\Dl \tau}-A_\tau
 U(\tau,t)\,x \\
 = A_{\tau-\theta\Dl\tau}U(\tau-\theta\Dl\tau,t)\,x
 - A_\tau U(\tau,t)\,x \to 0, \quad \Dl\tau\to 0,
 }
 where the right hand side converges to zero uniformly in $\tau\in[s,t]$
 since the continuous function
 $[s,t]\to E, \; \zeta \mto A_\zeta U(\zeta,t)$ is uniformly
 continuous.

 Thus, we have proved that $Q_{t_0,t_1}\ldots Q_{t_{n-1},t_n}x
 \to U(s,t)x$ as $n\to \infty$ for each $x\in Y$
 where $Y$ is dense in $E$.
 Note that by Assumption
 \ref{Qcond} of the theorem, the bounded on $E$ operators $Q_{t_0,t_1}\ldots
 Q_{t_{n-1},t_n}$ are bounded uniformly in
 $\{t_0, t_1, \ldots , t_n\}$.
 Hence, for an arbitrary subinterval $[s,t]\sub (S,T]$,
 the convergence $Q_{t_0,t_1}\ldots Q_{t_{n-1},t_n}x
 \to U(s,t)x$ holds
 for all $x\in E$.
 Thus, we have proved the theorem for the case $s>S$.

 Now we consider the case $s=S$. Let $s_N$ be a decreasing system of real numbers
 such that $\lim_{N\to\infty} s_N = s$. For each fixed $N$ and for all $x\in E$, we have:
 \eee{
  \label{08}
   Q_{t_0^N,t_1}\ldots Q_{t_{n-1},t_n} x \to U(s_N,t)\,x  \quad \text{as} \; \; n\to \infty
  }
  where $s_N=t_0^N < t_1 < \cdots < t_n = t$ is a partition of $[s_N,t]$.
  Note that for each fixed $\tau$ and for each fixed $x$,
  $\frac1{\Dl\tau} \|(Q_{\tau-\Dl\tau,\tau}-I)x\|_\e$ is bounded,
  which follows from convergence \eqref{ass_convergence} in Assumption \ref{ch_ass} if we set $t=\tau$.
  By the Banach-Steinhaus theorem there exists a constant $M_\tau > 0$ such that
  $\frac1{\Dl\tau} \|Q_{\tau-\Dl\tau,\tau}-I\|_{\sss E\to E} < M_\tau$. This implies
  that $\|Q_{\tau-\Dl\tau,\tau}-I\|_{\sss E\to E}$ tends to zero as $\Dl\tau\to 0$.
  Let us fix an arbitrary $\eps>0$, and find a $\dl > 0$ such that
  $\|Q_{s,s_N}-I\|_{\sss E\to E} < \eps$ and $\|U(s,s_N)-I\|_{\sss E\to E}<\eps$  whenever $s_N-s<\dl$.
  By Assumption \ref{Qcond} of Theorem \ref{chernnon},
  \mm{
  \|Q_{s,s_N}Q_{s_N,t_1} \ldots Q_{t_{n-1},t_n} - Q_{s_N,t_1} \ldots
  Q_{t_{n-1},t_n}\|_{\sss E\to E} \\
  \lt \|Q_{s,s_N}-I\|_{\sss E\to E} \,
  \|Q_{s_N,t_1} \ldots Q_{t_{n-1},t_n}\|_{\sss E\to E}
  \lt M(s,t) \, \eps.
  }
   By continuity of $U(\fdot,t)$,
  \mm{
  \|U(s,t)-U(s_N,t)\|_{\sss E\to E} \lt \|U(s,s_N)-I\|_{\sss E\to
  E}\, \|U(s_N,t)\|_{\sss E\to E}\\
 \lt
  \sup_{\xi\in [S,t]} \|U(\xi,t)\|_{\sss E\to E}\, \eps
  }
  Convergence \eqref{08} and two last estimates imply that for all $x\in E$
 \[
 Q_{t_0,t_1}\ldots Q_{t_{n-1},t_n} x \to U(s,t)\,x  \quad \text{as} \; \; n\to \infty.
 \]
The theorem is proved.
 \end{proof}
 \begin{lem}
 \label{lem5}
 Let $Y\sub \cap_{\tau\in [S,T]} D(A_\tau)$
 be a  Banach space, dense in $E$, let $[s,t]\sub (S,T]$ be a closed interval.
 Further
 let $B_{\tau,\Dl\tau}: Y\to E$,
 $\tau\in [s,t]$, $\Dl\tau \in (0,s-S)$,
 be bounded operators, and let
 $U(\xi,\tau)$, $(\xi,\tau)\in D_{[S,T]}$, be an evolution family of operators
 with the generators $A_\tau$. We assume that the following
 assumptions are fulfilled:
 \begin{enumerate}
 \item
 \label{33bounded}
 for every $y\in Y$,
 $\|B_{\tau,\Dl\tau}y\|_{\sss E}$ is bounded uniformly
 in $\tau\in [s,t]$ and
 $\Dl\tau \in [\dl,s-S]$,  where $\dl\in (0,s-S)$
 is fixed arbitrary;
 \item
 \label{A2}
 $U(\tau,t)Y\sub Y$ for all $\tau\in [s,t]$;
 \item
 \label{A3}
 for each $y\in Y$,
 the mapping $[s,t]\to Y, \; \tau\mto U(\tau,t)y$, is continuous;
\item
 \label{33}
 for each fixed $y\in Y$,
 \[
   \sup_{\tau\in [s,t]}\|A_\tau y\|_\e<\infty;
 \]
 \item
 \label{A333}
 for each fixed $y\in Y$,
 \eee{
 \label{333}
 \lim_{\Dl\tau\to 0} B_{\tau,\Dl\tau}y = A_\tau y
 }
 where the convergence is uniform in $\tau\in [s,t]$.
 \end{enumerate}
 Then, for every $y\in Y$,
 \eee{
  \label{3333}
  \lim_{\Dl\tau\to 0} B_{\tau,\Dl\tau}U(\tau,t)y = A_\tau
  U(\tau,t)y
 }
 and the convergence is uniform in $\tau\in [s,t]$.
 \end{lem}
 \begin{proof}
 Assumptions \ref{33} and \ref{A333}, along with Assumption
 \ref{33bounded}, imply that for each fixed $y\in Y$,
 $\|B_{\tau,\Dl\tau}y\|_{\sss E}$ is bounded uniformly in
 $\tau\in [s,t]$, and $\Dl\tau \in (0,s-S]$.
 By the Banach-Steinhaus theorem,
 $\|B_{\tau,\Dl\tau}\|_{\sss Y\to E}$
 are bounded uniformly in $\tau\in
 [s,t]$ and $\Dl\tau\in (0,s-S]$,
 i.e. there exists a constant $K$ such that
 \[
 \|B_{\tau,\Dl\tau}\|_{\sss Y\to E} < K.
 \]
 We fix a $y\in Y$. The set
 \eee{
 \label{comp3}
 \{U(\tau,t)y, \; \tau\in [s,t]\}
 }
 is a compact in $Y$ due to the continuity of the mapping
 $[s,t]\to Y, \; \tau\mto U(\tau,t)y$.
 Next, we fix an arbitrary small $\eps>0$ and
 find a finite $\eps$-net $\{y_i\}_{i=1}^N\sub Y$ for
 the compact \eqref{comp3}.
 Furthermore, we find a small $\dl>0$, such that
 for all $\tau\in [s,t]$, for all
 $\Dl\tau \in (0,\dl)$, and
 for all $y_i$, $1\lt i\lt N$,
 \[
 \|B_{\tau,\Dl\tau}y_i - A_\tau y_i\|_\e<\eps.
 \]
 Let $\tau\in [s,t]$ be fixed arbitrary, and $y_i$ be such that
 $\|U(\tau,t)y-y_i\|_{\e}<\eps$.
 We obtain:
 \[
 \|B_{\tau,\Dl\tau}U(\tau,t)y - B_{\tau,\Dl\tau}y_i\|_\e< K \|U(\tau,t)y-y_i\|_{\e}
 < K\eps.
 \]
 Taking the limit in the right hand side, as $\Dl\tau\to 0$,
 we obtain
 \[
 \|A_\tau U(\tau,t)y-A_\tau y_i\|_\e \lt K\eps.
 \]
 This implies:
 \aa{
 &\|B_{\tau,\Dl\tau}U(\tau,t)y - A_\tau U(\tau,t)y\|_\e \\ \lt
 &\|B_{\tau,\Dl\tau}U(\tau,t)y - B_{\tau,\Dl\tau}y_i\|_\e
 + \|B_{\tau,\Dl\tau}y_i - A_\tau y_i\|_\e
 + \|A_\tau y_i - A_\tau U(\tau,t)y\|_\e \\ <
 & (2K+1)\eps.
 }
 This proves that the limit \eqref{3333} exists
 and is uniform in $\tau\in[s,t]$. The lemma is proved.
 \end{proof}
 \begin{thm}[Corollary of Chernoff's Theorem]
 \label{chernoff_corollary}
 Let $A_t$, $Q_{t_1,t_2}$, $U(s,t)$, and $Y$ be as in
 Theorem \ref{chernnon}. Let us assume that Assumptions
 \ref{well-posedness}--\ref{Qcond} of Theorem \ref{chernnon}
 are fulfilled, and that
 for any subinterval $[s,t]\sub (S,T]$,
 Assumptions \ref{A3} and \ref{33} of Lemma \ref{lem5}
 are fulfilled. We assume, that for any $y\in Y$,
 \[
 \lim_{\Dl\tau\to 0}
 \frac{Q_{\tau-\Dl\tau,\tau}-I}{\Dl\tau} y = A_\tau y
 \]
 where the convergence is uniform in $\tau$
 running over closed subintervals $[s,t]\sub (S,T]$.
 Then, the statement of Theorem \ref{chernnon} holds
 true, i.e.
 for any subinterval $[s,t]\sub [S,T]$, for any sequence of
 partitions
 $\{s=t_0 < t_1 < \cdots < t_n=t\}$
 of $[s,t]$ such that
 $\max \,(t_{j+1}-t_j)\to 0$ as $n\to\infty$,
 and for all $x\in E$,
\[
 Q_{t_0,t_1}\ldots Q_{t_{n-1},t_n}x
 \to U(s,t)\,x, \quad n\to \infty.
 \]
 \end{thm}
 \begin{proof}
 Since we assume that Assumptions \ref{well-posedness}--\ref{Qcond} of Theorem
 \ref{chernnon} are  fulfilled, it suffices to prove that Assumption
 \ref{ch_ass} of Theorem \ref{chernnon} is fulfilled. This
 will follow from Lemma \ref{lem5} if we prove that
 Assumptions \ref{33bounded} and \ref{A333} of this lemma
 are fulfilled for
 the operators
 $B_{\tau,\Dl\tau}=\frac{Q_{\tau-\Dl\tau,\tau}-I}{\Dl\tau}$.
 Assumptions \ref{A2}--\ref{33} of Lemma \ref{lem5} clearly follow from
  those assumptions of Theorem \ref{chernnon} and  Lemma \ref{lem5}
 that are assumed here to be fulfilled.
 To prove Assumption \ref{ch_ass} of Theorem \ref{chernnon},  we fix an arbitrary closed interval
 $[s,t]\sub (S,T]$, and a
 $\dl\in (0,s-S)$. Then,
 for $\Dl\tau\in [\dl, s-S]$, we obtain:
 \[
 \Bigl\| \frac{Q_{\tau-\Dl\tau,\tau}-I}{\Dl\tau}\Bigr\|_{\sss E\to
 E} < \frac{M(s,t)+1}{\dl}
 \]
 where $M(s,t)$ is the constant in Assumption \ref{Qcond}
 of Theorem \ref{chernnon}. Assumption \ref{A333}  of Lemma
 \ref{lem5} is obviously fulfilled. By Lemma \ref{lem5},
 \[
 \lim_{\Dl\tau\to 0}\frac{Q_{\tau-\Dl\tau,\tau}-I}{\Dl\tau}U(\tau,t)y
 = A_\tau U(\tau,t)y
 \]
 and the limit is uniform in $\tau\in [s,t]$.
 Applying Theorem \ref{chernnon} completes the proof of the
 theorem.
 \end{proof}

 \subsection{Case of commuting generators}
 The following result has been obtained in \cite{sch}
 (p. 489, Proposition 2.5):
 \begin{pro}
 \label{pro1}
 Let $\{A_t\}$ be a stable family of pairwise commuting generators
 of strongly continuous semigroups. Let us assume that there exists a
 space $Y\sub \cap_{t\in [S,T]}D(A_t)$
 which is dense in $E$, and let for all $y\in Y$,
 the mapping $[S,T]\to E, \, t\mto A_ty$ be continuous. Then,
 $(\int_s^t A_r dr, Y)$
 is closable and its closure (which we
 still denote by $\int_s^t A_r dr$) is a generator. Moreover,
 the the Cauchy problem \eqref{Cauchy2} is well-posed and the
 evolution family solving \eqref{Cauchy2} is given by
 \[
 U(s,t)= e^{\int_s^t A_r dr}, \quad s\lt t\,.
 \]
 \end{pro}

 \begin{thm}[Chernoff's theorem for evolution families]
 \label{chern}
 Let $A_t$ be a stable family of pairwise commuting generators of
 strongly continuous semigroups, and let
 $Q_{t_1,t_2}$, $t_1,t_2>0$, be a
 two-parameter family of bounded operators $E\to E$,
 such that Assumptions \ref{continuous}--\ref{ch_ass} of Theorem \ref{chernnon}
 are fulfilled.
 Then, for any subinterval $[s,t]\sub [S,T]$, for any sequence of
 partitions
 $\{s=t_0 < t_1 < \cdots < t_n=t\}$
 of $[s,t]$ such that
 $\max \,(t_{j+1}-t_j)\to 0$ as $n\to\infty$,
 and for all $x\in E$,
 \[
 Q_{t_0,t_1}\ldots Q_{t_{n-1},t_n}x
 \to e^{\int_s^t A_r dr}\,x, \quad n\to \infty.
 \]
 \end{thm}
 \begin{proof}
 Proposition \ref{pro1} implies that Cauchy problem
 \eqref{Cauchy2} is well-posed, and that $U(s,t)= e^{\int_s^t A_r dr}$
 is the evolution family solving the Cauchy problem
 \eqref{Cauchy2}.
 Now the statement of the theorem follows immediately from Theorem~\ref{chernnon}.
 \end{proof}

 \section{Chernoff's theorem for evolution families
 generated by manifold valued stochastic processes}
 Let $M$ be a $\C^k$-smooth compact manifold,
 and let
 $A_0(t,x)$, $A_1(t,x)$, $\dots$, $A_d(t,x)$,
 $t\in [S,T]$, $x\in M$,
 be $\C^k$-smooth vector
 fields on $M$. This means that
 if $f\in \C^j(M)$ and $j>k$, then $A_i(t,\fdot)f\in C^k(M)$
 and if $j\lt k$, then
 $A_i(t,\fdot)f\in C^{j-1}(M)$ for all $t\in [S,T]$.
 Let us consider $t$-dependent second order differential operators:
 \eee{
 \label{gen_gen}
 A_t = \frac12 \sum_{\al=1}^d A_\al(t,\fdot) \circ
 A_\al(t,\fdot) + A_0(t,\fdot)
 }
 with the common domain $C^k(M)$
 independent of $t$.
 In the space $C^k(M)$ we introduce the norm:
 \eee{
 \label{k-norm}
\|f\|_{k} = \sum_{|\la|=0}^{k} \sup_y \sup_{x\in V}
 |\pl^\la f \circ \psi_y(x)|
 }
 where $\{(V,\psi_y), y\in M\}$ is an atlas covering
 $M$. The fact that $\|\fdot\|_{k}$ defines a norm is proved
 in \cite{Lawson} (pp. 175-176).
 The space $C^k(M)$ with the norm $\|\fdot\|_k$
 becomes a Banach space.
 We denote it by $Y$.

 Given a probability space $(\Om, \mc F, \mathbb P)$ with
 the filtration $\mc F_t$,
 and a $d$-dimensional $\mc F_t$-Brownian motion $B^\al_t$,
 we consider the stochastic differential equation:
 \aaa{
 \label{SDE}
 \begin{cases}
 & dX_t = A_\al(t,X_t)\circ dB^\al_t + A_0(t,X_t)dt\\
 & X_s = x
 \end{cases}
 }
 where $A_\al(t,X_t)\circ dB^\al_t$ is the Stratonovich
 differential. We denote by $\mathbb E$ the expectation
 relative to the measure $\mathbb P$.
 The operators $A_t$ are generators of diffusions $X_t$
 on $M$.
 \begin{lem}
  \label{lem8}
  Let $Y=(C^k(M), \|\fdot\|_k)$ where $k\gt 3$.
  Then, the solution of Cauchy problem \eqref{Cauchy2}
  on the interval $[S,t]$
  with the generators \eqref{gen_gen} and
  with the final condition  $u(t,x)=f(x)$,
  $f\in Y$, $x\in M$, exists, it is unique, and is
  given by
  \[
   u(s,x)= \mathbb E [f(X_t(s,x))]
  \]
  where $X_t(s,x)$
  is the solution of SDE \eqref{SDE}.
  Moreover, $u(s,x)\in Y$, $\frac{\pl}{\pl s} u(s,x) \in Y$ for all $s\in
  [S,t]$, and
  the mappings $[S,t]\to Y, s\mto u(s,\fdot)$, and
  $[S,t]\to E, s\mto \frac{\pl}{\pl s} u(s,\fdot)$
  are continuous.
  \end{lem}
  \begin{proof}
  Theorem 1.3 of Chapter 5 in \cite{Dalecky_Belopolskaya} (p. 182)
  implies as a particular case that
  \eee{
  \label{44}
   u(s,x)= (U(s,t)f)(x)=\mathbb E [f(X_t(s,x))]
  }
  is a solution of the Cauchy problem \eqref{Cauchy2}. Here
  $U(s,t)$,  $(s,t)\in D_{[S,T]}$, is the evolution family solving
  this Cauchy problem. Consider another evolution family $\td U(\tau,
  \xi)$, $(\tau,\xi)\in D_{[S,T]}$, satisfying the relation  $U(s,t)=\td
  U(T+S-t,T+S-s)$. Evidently, there exists another SDE of type \eqref{SDE} with
  $\C_k$- smooth coefficients, having a unique solution $\td
  X_\xi(\tau,x)$, such that for all $f\in Y$, for all $x\in M$,
  \[
  (\td U(\tau,\xi)f)(x)= \mathbb E [f(\td X_\xi(\tau,x))]\,.
  \]
  Applying Ito's formula gives:
  \mmm{
  \label{3-1}
  (U(s,t)f)(x) =( \td U (S+T-t,S+T-s)f)(x) = \mathbb E [f(\td
  X_{S+T-s}(S+T-t,x))] \\
  = f(x)+\int_0^{S+T-s} \mathbb E[(A_{S+T-\zeta}f)(\td
  X_\zeta(S+T-t,x))]d\zeta
  }
  where we have exchanged the symbol $\mathbb E$ for expectation
  with the integral in $\zeta$ by Fubini's theorem.
  For the partial derivative in $s$ we obtain:
  \eee{
  \label{3-2}
  \frac{\pl}{\pl s} u(s,x) = \frac{\pl}{\pl s} (U(s,t)f)(x) =  - \mathbb E[(A_{S+T-s}f)(\td
  X_s(S+T-t,x))]\,.
  }
 Clearly, $u(s,\fdot)\in Y$ and  $\frac{\pl}{\pl s}u(s,\fdot)\in
  Y$.
  Also, relations \eqref{3-1} and \eqref{3-2} imply that the
  mappings $[S,t]\to Y, s\mto u(s,\fdot)$ and
  $[S,t]\to E, s\mto \frac{\pl}{\pl s} u(s,\fdot)$
  are continuous. The lemma is proved.
  \end{proof}
 \begin{thm}[Chernoff's theorem for evolution
 families generated by manifold valued stochastic processes]
 \label{chern_application}
 Let $A_t$, $t\in [S,T]$, be given by \eqref{gen_gen},
 and let $D(A_t)=Y$ for all $t$.
 Further, let
 $Q_{t_1,t_2}$, $S\lt t_1 < t_2 \lt T$,
 be a family of contractions on $\C(M)$.
 We assume that the following assumptions are fulfilled:
 \begin{enumerate}
 \item
 \label{continuous_a}
 the functions $[S,T]\to \C(M)\,,\, t\mto A_tf$
 are continuously differentiable for all $f\in Y$;
 \item
 \label{AC2}
 stochastic differential equation \eqref{SDE}
 has a unique solution $X_t(s,x);$\footnote{Sufficient conditions under which
 \eqref{SDE} has a unique solution can be found for example in
 \cite{Dalecky_Belopolskaya} and \cite{Wat}}
 \item
 \label{ch_ass_a}
 for all $f\in Y$,
 \[
 \lim_{\Dl\tau\to 0}\frac{Q_{\tau-\Dl\tau,\tau}-I}{\Dl\tau}f
 = A_\tau f
 \]
 and the limit is uniform in $\tau$ running over closed intervals
 $[s,t]\sub (S,T]$.
 \end{enumerate}
 Then, for any subinterval $[s,t]\sub [S,T]$, for any sequence of
 partitions
 $\{s=t_0 < t_1 < \cdots < t_n=t\}$
 of $[s,t]$ such that
 $\max \,(t_{j+1}-t_j)\to 0$ as $n\to\infty$,
 and for all $f\in \C(M)$,
 the following convergence holds in $\C(M)$:
\[
 (Q_{t_0,t_1}\ldots Q_{t_{n-1},t_n}f)(\fdot)
 \to  \mathbb E[f(X_t(s,\fdot))], \quad n\to \infty.
 \]
  \end{thm}
  \begin{proof}
  Let $[s,t]\sub (S,T]$ be fixed.
  We would like to apply Theorem \ref{chernoff_corollary}.
  To this end, we have to verify Assumptions \ref{well-posedness} --
  \ref{Qcond} of Theorem \ref{chernnon} and Assumptions
  \ref{A2} and \ref{A3} of Lemma \ref{lem5}.
  Assumption \ref{well-posedness} of Theorem \ref{chernnon}
  follows from
  the paper \cite{kat53} by Kato.
  The paper \cite{kat53}
  guaranties existence and uniqueness of the solution
  of the Cauchy problem \eqref{Cauchy2} if the following assumptions are
  fulfilled:
  1) $D(A_t)=Y$ for all $t\in [S,T]$, and $Y$ is dense in $E$;
  2) the functions $t\mto A_t f$ are
  continuously differentiable.
  Due to this result,
  Assumption
  \ref{well-posedness} of Theorem \ref{chernnon} is fulfilled.
  Let $U(s,t)$ be the evolution family solving the Cauchy problem \eqref{Cauchy2},
  and let $u(s,x)$ denote the solution of  \eqref{Cauchy2}  with the final
  condition $f(x)$ at time $t$.
  Assumption \ref{continuous} of Theorem \ref{chernnon}
  is fulfilled by Lemma \ref{lem8}. Assumption \ref{Qcond}
  of Theorem \ref{chernnon} is fulfilled since
  $Q_{t_1,t_2}$ are contractions.
  Assumptions
  \ref{A2} and \ref{A3} of Lemma \ref{lem5}
  follow immediately from Lemma~\ref{lem8}.
  Now the statement of the theorem is implied by Theorem
  \ref{chernoff_corollary}.
  \end{proof}
  \section{Example: a time-inhomogeneous
  manifold valued stochastic process constructed by Chernoff's theorem}

 Below, we describe a construction of a
 time-inhomogeneous Markov process on a  compact Riemannian
 manifold using Theorem~\ref{chern_application}.
 Let $M$ be a compact Riemannian manifold without boundary
 isometrically imbedded into $\Rnu^m$, and  $\dim M=d$.
 Let $B_t$ be a Brownian motion on $\Rnu^m$ starting at the origin,
 and let
 $\ffi: [0,1] \to M$ be a two times continuously differentiable (non-random)
 function such that $\ffi(0)=x$.
 We consider the process $W_t = B_t + \ffi(t)$.
 Let $\W_\ffi$ be its distribution, $P_\ffi(t_1,z,t_2,A)$
 be its transition probability. Clearly,
 \aaa{
 \notag
 P_\ffi(t_1,z,t_2,A)= & P_0(t_1,z-\ffi(t_1), t_2, A-\ffi(t_2))\\
  \label{Pffi}
 =& \frac1{2\pi(t_2-t_1)^{\frac{m}2}}
 \int_A e^{-\frac{-|z-y-(\ffi(t_1)-\ffi(t_2))|^2}{2(t_2-t_1)}} dy
 }
 where $P_0$ corresponds to the case when $\ffi$
 is equal to zero identically.

 Let $U_\eps(M)$ be the $\eps$-neighborhood of $M$, and
 let $\W^x_{\eps,t}$ be the distribution of the process
 which is conditioned to take a value in $U_\eps(M)$ at time $t$.
 Specifically, we define a measure $\W^x_{\eps,t}$ by the
 following expression on the right hand side:
 \eee{
 \label{17}
 \int_{\C([0,t],\Rnu^m)} f(\om) \, \W^x_{\eps,t}(d\om)
 = \frac{\int_{\C([0,t],\Rnu^m)}\ind_{\{\om:\,\om(t)\in
 U_\eps(M)\}}(\om) \, f(\om) \,
   \W^x_{\ffi}(d\om)}
  {\W^x_{\ffi}\{\om:\,\om(t)\in
 U_\eps(M)\}}.
 }
  Let $P^x_{\eps,t}(\fdot,\fdot,\fdot,\fdot)$ be the transition probability for
  the distribution $\W^x_{\eps,t}$.
  By \eqref{Pffi} and \eqref{17}, $P^x_{\eps,t}(\fdot,\fdot,t,\fdot)$
  is given by
 \aaa{
 \label{10}
 \begin{split}
 \int_{\Rnu^m} P^x_{\eps,t}(s,z,t, dy) \, g(y) = &
 \int_{\C([0,t],\Rnu^m)} g(\om(t)) \W^x_{\eps,t}(d\om) \\
 = &\frac{
 \int_{U_\eps(M)} e^{-\frac{-|z-y-(\ffi(s)-\ffi(t))|^2}
 {2(t-s)}} g(y)\, dy
 }
 {
 \int_{U_\eps(M)} e^{-\frac{-|z-\bar y-(\ffi(s)-\ffi(t))|^2}
 {2(t-s)}}  d\bar y
 }
 \end{split}
 }
 where  $g:\Rnu^m\to\Rnu$ is bounded and continuous.
 Obviously, as $\eps\to 0$, the limit of the right hand side
 exists. Hence, the weak limit
 $P_{[s,t]}$ of the measures $P^x_{\eps,t}(s,\,\cdot\,,t,\,\cdot\,)$
 exists and equals
 \aa{
 \int_{\Rnu^m} P_{[s,t]}(z, dy) \, g(y) = &
 \frac{
 \int_M e^{-\frac{|z-y-(\ffi(s)-\ffi(t))|^2}
 {2(t-s)}} g(y)\, \la_M(dy)
 }
 {
 \int_M e^{-\frac{|z-\bar y-(\ffi(s)-\ffi(t))|^2}
 {2(t-s)}}  \la_M(d\bar y)
 } \\
 = & \int_M q_{\ffi}(s,z,t,y) g(y)\, \la_M(dy)
 }
 where $\la_M$ is the volume measure on $M$, and
 \[
 q_{\ffi}(s,z,t,y) =  \frac{
 e^{-\frac{|z-y-(\ffi(s)-\ffi(t))|^2}
 {2(t-s)}}
 }
 {
 \int_M e^{-\frac{|z-\bar y-(\ffi(s)-\ffi(t))|^2}
 {2(t-s)}}  \la_M(d\bar y)
 }\,.
 \]
 Given an interval $[s,t]$, the family of functions
 \eee{
 p_{\ffi}(t_1,z,t_2,y)=
 \frac1{2\pi(t_2-t_1)^{\frac{m}2}}
 \, e^{-\frac{|z-y-(\ffi(t_1)-\ffi(t_2))|^2}{2(t_2-t_1)}}
 }
 $s<t_1<t_2<t$, together with the function $q_{\ffi}(t_3,z,t,y)$, $t_3<t$,
 builds a family of transition densities that defines the distribution
 of a Markov process on $[s,t]$ conditioned
 to take a value on $M$ at time $t$.

 Consider a partition
 $\mc P = \{s=t_0 < t_1 < \cdots <t_n=t\}$. For each partition
 interval $[t_i,t_{i+1}]$, for each pair of points $\xi$ and $\tau$ such that
 $t_i<\xi<\tau\lt t_{i+1}$, and for each Borel set $A\sub\Rnu^m$, we
 define
 \eee{
 \label{Qdef}
Q(\xi,z,\tau,A)=
 \begin{cases}
 \int_A p_{\ffi}(\xi,z,\tau,y)\, dy, \quad  \tau<t_{i+1}, \\
 \int_{A\cap M} q_{\ffi} (\xi,z,\tau,y)\, \la_M(dy), \quad
 \tau=t_{i+1}.
 \end{cases}
 }
 Next, we
 add more points to the partition $\mc P$ to
 obtain a partition
 $\mc P'=\{s=\xi_0 < \xi_1< \cdots < \xi_N=t\}$
 containing $\mc P$.
 The family of measures
 \mm{
 Q(s,x,t,A) = \int_{\Rnu^m} Q(s,x,\xi_1,dx_1)
 \int_{\Rnu^m}Q(\xi_1,x_1,\xi_2,dx_2)
 \dots \\
 \int_{\Rnu^m}
 Q(\xi_{N-2},x_{N-2},\xi_{N-1},dx_{N-1})
 Q(\xi_{N-1},x_{N-1},t,A)
 }
 is a family of transition probabilities for a
 Markov process starting at the point $x\in M$
 at time $s$, and
 conditioned to take values on $M$ at
 all points of the partition $\mc P$.

  We apply Theorem \ref{chern_application} to a subfamily of the family
 $Q(\,\cdot\,, \,\cdot\,,\,\cdot\,, \,\cdot\,)$. Specifically, we
 investigate weak convergence of the family
\mm{
 q_{\sss \mc P}(s,x,t,y)=
 \int_M q_{\ffi}(s,x,t_1,x_1)\la_M(dx_1)
 \int_M q_{\ffi}(t_1,x_1,t_2,x_2)\la_M(dx_2)
 \dots \\
 \int_M q_{\ffi}(t_{n-2},x_{n-2},t_{n-1},x_{n-1})
  q_{\ffi}(t_{n-1},x_{n-1},t,y)\, \la_M(dx_{n-1})\, .
 }
 This family is a subfamily of $Q(\,\cdot\,, \,\cdot\,,\,\cdot\,, \,\cdot\,)$
 by definition \eqref{Qdef} of the family $Q$.
 We consider the following two
 parameter family of contractions $\C(M)\to\C(M)$:
 \eee{
  \label{Qq}
  (Q_{t_i,t_{i+1}}f)(\fdot)=\int_{M} q_\ffi (t_i,\fdot,t_{i+1},y) f(y)\la_M(dy).
  }
 \begin{thm}
 \label{application}
 As the mesh of $\mc P$ tends to zero,
 the following convergence holds in $\C(M)$:
 \eee{
 \label{01}
 \int_M q_{\sss \mc P}(s,\fdot,t,y)\, g(y)\, \la_M(dy)
 \to \int_M p(s,\fdot,t,y)\, g(y) \, \la_M(dy)
 }
 where  $g\in \C(M)$, $p(s,x,t,y)$
 is the transition density function of the
 process generated by
 \eee{
 \label{generator}
 A_s = (\ffi'(s),\nab_M))_{\Rnu^m}-\frac12\lap_M.
 }
 \end{thm}
 \begin{lem}
 The $A_s$ given by \eqref{generator} generate contraction
 semigroups on $\C(M)$. Moreover, each $A_s$ is the generator
 of a diffusion $X(\tau)$ on $M$ which is the solution of the following SDE:
 \eee{
 \label{st}
 \left\{
 \begin{split}
 & dr(\tau)=\td L_\al(r(\tau))\circ dw^\al_{s,\ffi},\\
 & r(0)=r,
 \end{split}
 \right.
 }
 where $r(\tau)=(X^i(\tau),e^i_\al(\tau))$, $\{e_\al(\tau)\}$ is a basis in the tangent
 space at the point $X(\tau)$,
 $\td L_\al$ are canonical horizontal vector fields~\cite{Wat},
 $w^\al_{s,\ffi}(\tau)=\ffi'(s)^\al \tau+B^\al(\tau)$, $B^\al(\tau)$
 is a Brownian motion in $\Rnu^d$.
 \end{lem}
 \begin{proof}
 Let $r(\tau)=(X^i(\tau),e_\al^i(\tau))$ be the solution of~\eqref{st}.
 We find the generator of $X(\tau)$.
 Consider the function $f(r)=f(x)$ for $r=(x,e)$.
 We have:
 \aa{
 f(X(\tau))-f(X(0))&=f(r(\tau))-f(r(0))\\
 &=\intl_0^\tau(\td L_\al f)(r(\xi))\circ dw_{s,\ffi}^\al(\xi)\\
 &=\intl_0^\tau\td L_\al f(r(\xi))dB^\al(\xi)
 +\intl_0^\tau\td L_\al f(r(\xi))\ffi'(s)^\al d\xi\\
 &\quad +\frac12\intl_0^\tau
 \suml_{\al=1}^d\td L_\al(\td L_\al f)(r(\xi))d\xi.
 }
 The definition of the generator of a process gives:
 \[
 A_s f=\suml_{\al=1}^d (\td L_\al f, \ffi'(s)^\al)_{\Rnu^d}
 +\frac12\suml_{\al=1}^d\td L_\al(\td L_\al f).
 \]
 Since $f(r)=f(x)$, i.e. does not depend on $e$, then
 the scalar product in the first term of the right hand side is
 well-defined,
 and
 \[
 \td L f = \nab_M f
 \]
 by definition of $\td L$.
 Further, it was shown in~\cite{Wat} (Chapter V, paragraph 3) that
 \[
 \suml_{\al=1}^d\td L_\al(\td L_\al f) = -\lap_M f\,.
 \]
 Thus, we have proved that
 \[
 A_s f=(\ffi'(s),\nab_M f)_{\Rnu^d}-\frac12\lap_M f = (\ffi'(s),\nab_M f)_{\Rnu^m}-\frac12\lap_M f\,.
 \]
 Since $A_s$ is  a generator of a diffusion on $M$, $A_s$ generates
 a contraction semigroup on $\C(M)$. The lemma is proved.
 \end{proof}
 \begin{proof}[Proof of Theorem \ref{application}]
 We apply Theorem \ref{chern_application} to the generators \eqref{generator}
 and the two-parameter family of contractions
 \[
 (Q_{t_1,t_2}f)(x) = \int_M q_{\ffi}(t_1,x,t_2,y)\, f(y)\,
 \la_M(dy).
 \]
 Assumption \ref{continuous_a} of Theorem \ref{chern_application}
 is fulfilled by continuity
 of $\ffi''(t)$ which we have assumed.
 Assumption \ref{AC2} of Theorem \ref{chern_application}
 is fulfilled by Theorem 2.1, p. 152, from the book
 \cite{Dalecky_Belopolskaya}, where the authors have considered a
 more general case of a manifold in a Banach space.
  We show that Assumption  \ref{ch_ass_a} of Theorem
  \ref{chern_application} is fulfilled too.
  We have to prove that
  \mmm{
  \label{2prove}
 \lim_{\Dl\tau\to 0}\frac1{\Dl\tau}\bigl
 (\int_M q_\ffi(\tau-\Dl\tau,y,\tau,z) \, g(z)\, \la_M(dz)  -
 g(y)\bigr)\\
 = (\ffi'(\tau),\nab_M g(y))_{\Rnu^m}-\frac12\, \lap_M g(y)\, ,
  }
 and that the limit is uniform in $\tau$.
 Introduce the notation: $\Dl\ffi_\tau= \ffi(\tau)-\ffi(\tau -\Dl\tau)$.
 \begin{lem}
 \label{lem6}
 Let  $g\in C^3(M)$.
 There exist a $\dl>0$, a constant $K_g>0$, and a function
 $R:[0,\dl]\x M\x C^3(M) \to \Rnu$ satisfying:
 \eee{
 \label{Rest}
  |R(\Dl\tau, \fdot, g)| < K_g \, \Dl\tau^{\frac12},
 }
 and such that for all $y\in M$
 the following relation holds:
 \aaa{
 \label{asymptotic_main}
 \begin{split}
 \frac{\int_M
 g(z)e^{-\frac{|z-y-\Dl\ffi_\tau|^2}{2\Dl\tau}}\la_M(dz)}
 {
  \int_M
  e^{-\frac{|z-y-\Dl\ffi_\tau|^2}{2\Dl\tau}}\la_M(dz)}
  = & \, g(y) + (\Dl\ffi_\tau,\nab_M g(y))_{\Rnu^m} - \frac{\Dl\tau}2\lap_Mg(y) \\
  + & \Dl\tau\, R(\Dl\tau,y, g).
  \end{split}
 }
 \end{lem}
 \begin{proof}
  We find a $U_\eps(M)$, the $\eps$-neighborhood of $M$,
 where
 the normal spaces
 $N_{y_1}$ and $N_{y_2}$ do not intersect each other for each
 pair of points $y_1\in M$ and $y_2\in M$.
 Hence,
 each $y\in U_\eps(M)$ can be uniquely presented as $y=z+tn(z)$,
 where $z\in M$, $n(z)\in N_z$, and $|n(z)|=1$.
 Let $P_M: U_\eps(M)\to M$, $z+tn(z)\mto z$,  $t\in \Rnu$,
 be the projection on $M$.
 For an arbitrary $u\in \Rnu^m$, $|u|< \eps$, and a $y\in M$, we define:
 \[
 u_M(y)= P_M(y+u)-u, \qquad
 u_{\bot}(y)=u-u_M(y).
 \]
 We have:
 \eee{
 \label{7}
  \frac{\int_M
  g(z)e^{-\frac{|z-y-\Dl\ffi_\tau|^2}{2\Dl\tau}}\la_M(dz)
  }
  {
  \int_M
  e^{-\frac{|z-y-\Dl\ffi_\tau|^2}{2\Dl\tau}}\la_M(dz)
  }
  =
  \frac{\int_M
  g(z)e^{-\frac{|z-y-(\Dl\ffi_\tau)_M(y)|^2}{2\Dl\tau}}
  e^{\frac{(z-y,(\Dl\ffi_\tau)_{\hspace{-0.7mm}\bot}(y))}{\Dl\tau}}
  \la_M(dz)
  }
  {
  \int_M
  e^{-\frac{|z-y-(\Dl\ffi_\tau)_M(y)|^2}{2\Dl\tau}}
  e^{\frac{(z-y,(\Dl\ffi_\tau)_{\hspace{-0.7mm}\bot}(y))}{\Dl\tau}}
  \la_M(dz)
  }\,.
 }
  We will need the following formula
 (see \cite{sm_vr} and \cite{mine4}):
 \eee{
 \label{sm_vr_fm}
 \frac{
 \int_M
 e^{-\frac{|z-y|^2}{2t}} \, h(z)\, \la_M(dz)
  }
  {
  \int_M
  e^{-\frac{|z-y|^2}{2t}}\la_M(dz)
  }
 = h(y)-\frac{t}2\,\lap_M h(y)  + t\bar R(t, y , h)
 }
 where $|\bar R(t,y,h)|< K\, \|h\|_3\,t^{\frac12}$, $K$ is a constant,
 and the norms $\|\fdot\|_3$ are defined by \eqref{k-norm} for $k=3$.
 Dividing the numerator and denominator in
 the right hand side of \eqref{7} by
 $\int_M e^{-\frac{|z-y-(\Dl\ffi_\tau)_M(y)|^2}{2\Dl\tau}}
 \la_M(dz)$,
 and applying \eqref{sm_vr_fm},
 we obtain:
  \mmm{
  \label{88}
  \frac{\int_M
  g(z)e^{-\frac{|z-y-\Dl\ffi_\tau|^2}{2\Dl\tau}}\la_M(dz)
  }
  {
  \int_M
  e^{-\frac{|z-y-\Dl\ffi_\tau|^2}{2\Dl\tau}}\la_M(dz)
  }
  \\
  =
  \frac{
   g\bigl(y+(\Dl\ffi_\tau)_M(y)\bigr) - \frac{\Dl\tau}2\,
   \lap_M \bigl( g\,
   e^{\frac{(\fdot -y,(\Dl\ffi_\tau)_{\hspace{-0.7mm}\bot}(y))}{\Dl\tau}}\bigr)(y)
   +\Dl\tau R_1
  }
  {1-\frac{\Dl\tau}2\, \bigl(\lap_M
  e^{\frac{(\fdot -y,(\Dl\ffi_\tau)_{\hspace{-0.7mm}\bot}(y))}{\Dl\tau}}
  \bigr)(y)
  + \Dl\tau R_2
  }
  }
  where $R_1$ and $R_2$ are short-hand notations for
  $\bar R(\Dl\tau,y+\Dl\ffi_\tau,g\,
  e^{\frac{(\fdot -y,(\Dl\ffi_\tau)_{\hspace{-0.7mm}\bot}(y))}{\Dl\tau}})$
  and
  $\bar R(\Dl\tau,y+\Dl\ffi_\tau,
  e^{\frac{(\fdot
  -y,(\Dl\ffi_\tau)_{\hspace{-0.7mm}\bot}(y))}{\Dl\tau}})$,
  respectively. They can be estimated as follows:
  \[
  |R_1| < \, \td K \|g\|_3
  \,\|
  e^{\frac{(\fdot -y,(\Dl\ffi_\tau)_{\hspace{-0.7mm}\bot}(y))}{\Dl\tau}}
  \|_3
  \, (\Dl\tau)^{\frac12},
   \quad
  |R_2| <  \, K \|
  e^{\frac{(\fdot -y,(\Dl\ffi_\tau)_{\hspace{-0.7mm}\bot}(y))}{\Dl\tau}}
  \|_3 \, (\Dl\tau)^{\frac12}.
  \]
  We show that $\bigl(\lap_M
  e^{\frac{(\fdot -y,(\Dl\ffi_\tau)_{\hspace{-0.7mm}\bot}(y))}{\Dl\tau}}
  \bigr)(y)$
  and
  $\|
  e^{\frac{(\fdot -y,(\Dl\ffi_\tau)_{\hspace{-0.7mm}\bot}(y))}{\Dl\tau}}
  \|_3$
  are bounded in $\tau$ and $\Dl\tau$.
  We consider local normal charts $\psi(\bar\xi)=\psi(\xi_1,\dots\xi_d)$
  at the point $y$.
  We obtain:
  \eee{
  \label{this}
  \frac{\pl}{\pl \xi_i}\,
  e^{\frac{(\psi(\bar\xi) -y,(\Dl\ffi_\tau)_{\hspace{-0.7mm}\bot}(y))}{\Dl\tau}}
  =
  e^{\frac{(\psi(\bar\xi)-y,
  (\Dl\ffi_\tau)_{\hspace{-0.7mm}\bot}(y))}{\Dl\tau}}\,
  \Bigl(\frac{\pl}{\pl \xi_i}\psi(\bar\xi),
  \frac{(\Dl\ffi_\tau)_{\hspace{-0.7mm}\bot}(y)}{\Dl\tau}
  \Bigr)_{\Rnu^m}.
  }
  This formula makes obvious the resulting expression
  upon taking two further derivatives, and thus, it
  shows that
  $\bigl(\lap_M
  e^{\frac{(\fdot -y,(\Dl\ffi_\tau)_{\hspace{-0.7mm}\bot}(y))}{\Dl\tau}}
  \bigr)(y)$
  and
  $\|
  e^{\frac{(\fdot -y,(\Dl\ffi_\tau)_{\hspace{-0.7mm}\bot}(y))}{\Dl\tau}}
  \|_3$
  are bounded in $\tau$ and $\Dl\tau$ if and only if
  $\frac{(\Dl\ffi_\tau)_{\hspace{-0.7mm}\bot}(y)}{\Dl\tau}$
  is bounded in $\tau$ and $\Dl\tau$.
  The latter fact holds by existence of the limit:
  \eee{
  \label{0}
  \lim_{\Dl\tau\to 0}\frac{(\Dl\ffi_\tau)_{\hspace{-0.7mm}\bot}(y)}{\Dl\tau}
  =\Pr\nolimits_{N_y}\ffi'(\tau)
  }
  where $\Pr\nolimits_{N_y}$ is the orthogonal projection onto $N_y$, the
  normal space at $y$. Now we can apply the short time asymptotic
  in $\Dl\tau$ to the denominator at the right hand side of
  \eqref{88}
  while using the relation
  $
  \lap_M \bigl( g\,
  e^{\frac{(\fdot -y,(\Dl\ffi_\tau)_{\hspace{-0.7mm}\bot}(y))}
  {\Dl\tau}}\bigr)
  = e^{\frac{(\fdot -y,(\Dl\ffi_\tau)_{\hspace{-0.7mm}\bot}(y))}
  {\Dl\tau}} \lap_M g
  -2\, (\nab_M
  e^{\frac{(\fdot -y,(\Dl\ffi_\tau)_{\hspace{-0.7mm}\bot}(y))}
  {\Dl\tau}}, \nab_M g)_{\Rnu^m}
  +g \, \lap_M
  e^{\frac{(\fdot -y,(\Dl\ffi_\tau)_{\hspace{-0.7mm}\bot}(y))}
  {\Dl\tau}}
  $.
  We obtain:
  \mmm{
  \label{00}
  \frac{\int_M
  g(z)e^{-\frac{|z-y-\Dl\ffi_\tau|^2}{2\Dl\tau}}\la_M(dz)
  }
  {
  \int_M
  e^{-\frac{|z-y-\Dl\ffi_\tau|^2}{2\Dl\tau}}\la_M(dz)
  }
  = g(y)+ \bigl((\Dl\ffi_\tau)_M(y),\nab_M g(y)\bigr)_{\Rnu^m}
  -\frac{\Dl\tau}2\, \lap_M g(y) \\
  +\Dl\tau \Bigl(\left.\nab_M
  e^{\frac{(\fdot -y,(\Dl\ffi_\tau)_{\hspace{-0.7mm}\bot}(y))}
  {\Dl\tau}}\right|_y, \nab_M g(y)\Bigr)_{\Rnu^m} + \Dl\tau \td R(\Dl\tau)
  }
  where $|\td R(\Dl\tau)| < \td K \|g\|_3 \Dl\tau^{\frac12}$,
  $\td K$ is a constant.
  Formulas \eqref{this} and \eqref{0} imply:
  \[
   \Bigl(
   \left.\nab_M\,
   e^{\frac{(\fdot -y,(\Dl\ffi_\tau)_{\hspace{-0.7mm}\bot}(y))}{\Dl\tau}}
   \right|_y,
   \nab_M g(y)
   \Bigr)
  =
  \Bigl(
  \frac{(\Dl\ffi_\tau)_{\hspace{-0.7mm}\bot}(y)}{\Dl\tau},
  \nab_M g(y)
  \Bigr).
  \]
  Substituting this in \eqref{00}, we obtain
  \eqref{asymptotic_main}.
 The lemma is proved.
 \end{proof}
 We continue the proof of Theorem \ref{application}. Lemma \ref{lem6} easily implies  the convergence in \eqref{2prove}.
 This convergence is uniform in $\tau$. Indeed, \eqref{Rest}
 implies that $R(\Dl\tau,\fdot,g)$ converges to zero uniformly
 in $\tau$. Further, we have:
 \[
 \Bigl|\frac{\Dl\ffi_\tau}{\Dl\tau} - \ffi'(\tau)\Bigr| \lt
 \frac12 \sup_{\theta\in
 [0,1]}|\ffi''(\tau+\theta\Dl\tau)|\Dl\tau,
 \]
 and hence, $\frac{\Dl\ffi_\tau}{\Dl\tau}$
 converges to $\ffi'(\tau)$ uniformly in $\tau$.
 Thus, we have verified all the assumptions
 of Theorem \ref{chern_application}.
 Finally, we note that since $p(\fdot,\fdot,\fdot,\fdot)$
 is the transition density for the diffusion process
 $X_t(s,x)$, and hence
 \eee{
 \label{Ep}
 \int_M p(s,x,t,y)f(y) \la_M(dy)
 = \mathbb E[f(X_t(s,x))].
 }
 Now the statement of the theorem is implied by convergence \eqref{01},
 formulas \eqref{Ep} and \eqref{Qq}, and by Theorem \ref{chern_application}.
 \end{proof}


\begin{thebibliography}{99}
 \bibitem{sm_vr}
 {\em Smolyanov O.G., Weizs\"acker H.v., Wittich O.}, Brownian motion
 on a manifold as a limit of stepwise conditioned standard
 Brownian motions, Canadian Mathematical Society, Conference Proceedings,
 Vol. 29, 2000, pp. 589-602.
 \bibitem{mine4}
 {\em \'E.\ Yu.\ Shamarova},
 Constructing a Brownian sheet with values in a compact Riemannian manifold, (Russian, English),
 Math. Notes 76, No. 4, 590-596 (2004); translation from Mat. Zametki 76, No. 4,
 635-640, 2004.
 \bibitem{chernoff}
 {\em Chernoff, R.}, Product formulas, Nonlinear semigroups,
 and Addition of Unbounded operators, Memoirsof American Math.
 Soc., 140, 1974
 \bibitem{Dalecky_Belopolskaya}
 {\em Dalecky Yu. L, Belopolskaya Ya. I.}, Stochastic equations
 and differential geometry, Series: Mathematics and its Applications,
 Kluwer Academic Publishers, Netherlands, 260 p., 1989
 \bibitem{Dalecky-Fomin}
 {\em Dalecky Yu. L, Fomin S. V.}, Measures and differential
 equations in infinite-dimensional space, Series: Mathematics and its Applications,
 Kluwer Academic Publishers, Netherlands, 337 p., 1991
 \bibitem{Egorov}
 {\em Egorov Yu. V., Shubin M.A.}, Partial Differential Equations,
 Vol III, Sprinter-Verlag
 \bibitem{nag}
 {\em Engel, K.-J., Nagel R.},
 One-parameter semigroups for linear evolution equations. (English)
 Graduate Texts in Mathematics. 194. Berlin: Springer. xxi, 586
 p., 2000.
 \bibitem{gold}
 {\em Goldstein, Jerome A.},
 Semigroups of linear operators and applications. (English) Oxford
 Mathematical Monographs. New York: Oxford University Press;
 Oxford: Clarendon Press. X, 245 p., 1985.
\bibitem{Lawson}
 Lawson, H.\ B., Michelsohn M.-L., Spin Geometry, Princeton University
 Press, Princeton NJ, 1989.
 \bibitem{tan}
 {\em Tanabe, Hiroki},
 Equations of evolution. Translated from Japanese by N. Mugibayashi
 and H. Haneda. (English) Monographs and Studies in Mathematics. 6.
 London - San Francisco - Melbourne: Pitman. XII, 1979, 260 p.
 \bibitem{sch}
 {\em Nickel, Gregor; Schnaubelt, Roland},
 An extension of Kato's stability condition for nonautonomous
 Cauchy problems. Taiwanese J. Math. 2, No.4, 483-496, 1998.
 \bibitem{kat53}
 {\em Kato, Tosio},
 Integration of the equation of evolution in a Banach space.
(English) J. Math. Soc. Japan 5, 208-234 (1953).
\bibitem{Wat}
{\em Ikeda N., Watanabe S.}, Stochastic differential equations and
 diffusion processes, North Holland publishing company (1989).
 \end{thebibliography}
\end{document}